\documentclass{article}
\usepackage{amssymb}
\usepackage{amsmath}
\usepackage{amstext}

\newtheorem{theorem}{Theorem}

\newtheorem{definition}[theorem]{Definition}

\newtheorem{remark}[theorem]{Remark}

\input tcilatex
\begin{document}

\title{A model problem for ultrafunctions}
\author{Vieri Benci\thanks{
Dipartimento di Matematica, Universit\`{a} degli Studi di Pisa, Via F.
Buonarroti 1/c, Pisa, ITALY and Department of Mathematics, College of
Science, King Saud University, Riyadh, 11451, SAUDI ARABIA. e-mail: \texttt{%
benci@dma.unipi.it}}, Lorenzo Luperi Baglini\thanks{
Dipartimento di Matematica, Universit\`{a} degli Studi di Pisa, Via F.
Buonarroti 1/c, Pisa, ITALY, e-mail: \texttt{lorenzo.luperi@for.unipi.it}}}
\maketitle
\date{}

\begin{abstract}
The purpose of this paper is to show that Non-Archimedean Mathematics (NAM),
namely mathematics which uses infinite and infinitesimal numbers,\ is useful
to model some Physical problems which cannot be described by the usual
mathematics. The problem which we will consider here is the minimization of
the functional%
\begin{equation*}
E(u,q)=\frac{1}{2}\int_{\Omega }|\nabla u(x)|^{2}dx+u(q).
\end{equation*}%
If $\Omega \subset \mathbb{R}^{N}$ is a bounded open set and $u\in \mathcal{C%
}_{0}^{2}(\overline{\Omega }),$ this problem has no solution since $\inf \
E(u,q)=-\infty .$ On the contrary, as we will show, this problem is well
posed in a suitable non-Archimedean frame. More precisely, we apply the
general ideas of NAM and some of the techniques of Non Standard Analysis to
a new notion of generalized functions, called ultrafunctions, which are a
particular class of functions based on a Non-Archimedean field. In this
class of functions, the above problem is well posed and it has a solution.

\medskip


\medskip

\end{abstract}

\tableofcontents

\section{Introduction}

By Non-Archimedean Mathematics (NAM) we mean mathematics based on
non-Archimedean fields, namely ordered fields which contain infinite and
infinitesimal numbers. We believe that, in many circumstances, NAM allows to
construct models of the physical world in a very elegant and simple way. In
the years around 1900, NAM was investigated by prominent mathematicians such
as Du Bois-Reymond \cite{DBR}, Veronese \cite{veronese}, David Hilbert \cite%
{hilb} and Tullio Levi-Civita \cite{LC}, but then it has been forgotten
until the '60s when Abraham Robinson presented his Non Standard Analysis
(NSA) \cite{rob}. We refer to Ehrlich \cite{el06} for a historical analysis
of these facts and to Keisler \cite{keisler76} for a very clear exposition
of NSA (see also \cite{benci99}, \cite{BDN2003})

The purpose of this paper is to show that NAM\ is useful to model some
Physical problems which cannot be described by the usual mathematics even if
they are relatively simple.

The notion of material point is a basic tool in Mathematical Physics since
the times of Euler who introduced it. Even if material points do not exist,
nevertheless they are very useful in the description of nature and they
semplify the models so that they can be treated by mathematical tools.
However, as new notions entered in Physics (such as the notion of field),
the use of material points led to situations which required new mathematics.
For example, in order to describe the electric field generated by a charged
point, we need the notion of Dirac measure $\delta _{q}$, namely this field
satifies the following equation:%
\begin{equation}
\Delta u=\delta _{q}  \label{one}
\end{equation}%
where $\Delta $ is the Laplace operator.

In this paper, we will describe a simple problem whose modelization requires
NAM. Let $\Omega \subseteq \mathbb{R}^{2}$ be an open bounded set which
represents a (ideal) membrane. Suppose that in $\Omega $ is placed a
material point $P$, which is left free to move.

Suppose that the point has a unit weight and the only forces acting on it
are the gravitational force and the reaction of the membrane. If $q\in
\Omega $ is the position of the point and $u(x)$ represents the profile of
the membrane, it follows that equation (\ref{one}) holds in $\Omega $ with
boundary condition $u=0$ on $\partial \Omega .$

The question is: which is the point $q_{0}\in \Omega $ that the particle
will occupy?\newline
The natural way to approach this problem would be the following: for every $%
q\in \Omega $, the energy of the system is given by the elastic energy plus
the gravitational energy, namely%
\begin{equation}
E(u,q)=\frac{1}{2}\int_{\Omega }|\nabla u(x)|^{2}dx+u(q)  \label{integ}
\end{equation}%
If the couple $(u_{0},q_{0})$ minimizes $E$, then $q_{0}$ is the equilibrium
point. For every $q\in \Omega $, let $u_{q}(x)$ be the configuration when $P$
is placed in $q,$ namely the solution of equation (\ref{one}). So the
equilibrium point $q_{0}$ is the the point in which the function

\begin{equation}
F(q)=E(u_{q},q)
\end{equation}%
has a minimum.

In the classical context, this "natural" approach cannot be applied; in fact 
$u_{q}(x)$ has a singularity at the point $q$ which makes $u(q)$ not well
defined and the integral in (\ref{integ}) to diverge. On the contrary, this
problem can be treated in NAM as we will show. In fact, since infinite
numbers are allowed, we will be able to find a minimum configuration for the
energy (\ref{integ}).

In order to pursue this program, we apply the general ideas of NAM and some
of the techniques of NSA to a new notion of generalized functions which we
have called \textbf{ultrafunctions} (see \cite{ultra}). Ultrafunctions are a
particular class of functions based on a superreal field $\mathbb{R}^{\ast
}\supset \mathbb{R}$. More exactly, to any continuous function $f:\mathbb{R}%
^{N}\rightarrow \mathbb{R}$ we associate in a canonical way an ultrafunction 
$f_{\Phi }:\left( \mathbb{R}^{\ast }\right) ^{N}\rightarrow \mathbb{R}^{\ast
}$ which extends $f;$ but the ultrafunctions are much more than the
functions and among them we can find solutions of functional equations such
as equation (\ref{one}) which are defined in every point of $\Omega ^{\ast
}\subset \left( \mathbb{R}^{\ast }\right) ^{N}.$ Thus, the energy (\ref%
{integ}) is well defined for every ultrafunction even if it might assume
infinite values.

Now we itemize some of the peculiar properties of the ultrafunctions:

\begin{itemize}
\item the space of ultrafunctions is larger than the space of distributions,
namely, to every distribution $T,$ we can associate in a canonical way an
ultrafunction $T_{\Phi }$ (for details see \cite{ultra}); in particular the
Dirac measure can be represented by an ultrafunction $\delta _{q}(x)$ and,
for every ultrafunction $u,$ we have that 
\begin{equation*}
\int u(x)\delta _{q}(x)dx=u(q);
\end{equation*}

\item similarly to the distributions, the ultrafunctions are motivated by
the need of having generalized solutions and also by the need to model
extreme physical situations which cannot be described by functions defined
in $\mathbb{R}^{N}$; however, while the distributions are no longer
functions, the ultrafunctions are still functions even if they have larger
domain and range;

\item unlike the distributions, the space of ultrafunctions is suitable for
non linear problems such as the one described above;

\item if a problem has a unique classical solution $u,$ then $u_{\Phi }$ is
the only solution in the space of ultrafunctions;

\item the main strategy to prove the existence of generalized solutions in
the space of ultrafunction is relatively simple; it is just a variant of the
Faedo-Galerkin method.
\end{itemize}

Before concluding the introduction, we refer to \cite{BGG} and to \cite{BHW}
where other situations which require NAM are presented.

\subsection{Notations}

Let $\Omega $\ be a subset of $\mathbb{R}^{N}$: then

\begin{itemize}
\item $\mathcal{C}\left( \Omega \right) $ denotes the set of real continuous
functions defined on $\Omega ;$

\item $\mathcal{C}_{0}\left( \overline{\Omega }\right) $ denotes the set of
real continuous functions on $\overline{\Omega }$ which vanish on $\partial
\Omega ;$

\item $\mathcal{C}^{k}\left( \Omega \right) $ denotes the set of functions
defined on $\Omega \subset \mathbb{R}^{N}$ which have continuous derivatives
up to the order $k;$

\item $\mathcal{C}_{0}^{k}\left( \overline{\Omega }\right) =\mathcal{C}%
^{k}\left( \overline{\Omega }\right) \cap \mathcal{C}_{0}\left( \overline{%
\Omega }\right) ;$

\item $\mathcal{D}\left( \Omega \right) $ denotes the set of the infinitely
differentiable functions with compact support defined on $\Omega \subset 
\mathbb{R}^{N};$

\item $H^{1}(\Omega )$ is the usual Sobolev space defined as the set of
functions $u\in L^{2}\left( \Omega \right) $ such that $\nabla u\in
L^{2}\left( \Omega \right) ;$

\item $H_{0}^{1}(\Omega )$ is the closure of $\mathcal{D}\left( \Omega
\right) $ in $H^{1}(\Omega );$

\item $H^{-1}(\Omega )$ is the topological dual of $H_{0}^{1}(\Omega ).$
\end{itemize}

\section{The ultrafunctions}

In this section we briefly recall the notion of $\Lambda $-limit and of
ultrafunction which have been introduced in \cite{ultra}.

\subsection{The $\Lambda$-limit}

The idea behind the concept of $\Lambda $-limit is the following: let $%
\mathbb{U}$ denote a "mathematical universe" (which will be precisely
introduced in definition (\ref{universe})), and $\mathcal{F}$ the set of
finite subsets of $\mathbb{U}$, ordered by inclusion. The $\Lambda $-limit
can be thought as a way to associate to every net $\varphi :\mathcal{F}%
\rightarrow \mathbb{R}$ a limit $\lim\limits_{\lambda \uparrow \mathbb{U}%
}\varphi (\lambda )$ that satisfies a few properties of coherence.\newline
These limits will be elements of a Non-Archimedean field $\mathbb{K}$; since
this leads to work in such fields, we recall a few basic facts and
definitions:

\begin{definition}
Let $\mathbb{K}$ be an ordered field. Let $\xi \in \mathbb{K}$. We say that:

\begin{itemize}
\item $\xi $ is infinitesimal if for all $n\in \mathbb{N}$ $|\xi |\ <\frac{1%
}{n}$;

\item $\xi $ is finite if there exists $n\in \mathbb{N}$ such as $|\xi |<n$;

\item $\xi $ is infinite if for all $n\in \mathbb{N}$ $|\xi |>n$.
\end{itemize}
\end{definition}

\begin{definition}
An ordered field $\mathbb{K}$ is called non-Archimedean if it contains an
infinitesimal $\xi \neq 0$.
\end{definition}

We are interested in fields that extend $\mathbb{R}$:

\begin{definition}
A superreal field is an ordered field $\mathbb{K}$ that properly extends $%
\mathbb{R}$.
\end{definition}

Since $\mathbb{R}$ is complete, it is easily seen that every superreal field
contains infinitesimal and infinite numbers.\newline
In order to precise the notion of $\Lambda $-limit, we need to define the
notion of "mathematical universe". For our applications, we take as
mathematical universe the superstructure on $\mathbb{R}$:

\begin{definition}
\label{universe} The \textbf{superstructure} on $\mathbb{R}$ is 
\begin{equation*}
\mathbb{U}=\dbigcup_{n=0}^{\infty }\mathbb{U}_{n}
\end{equation*}%
where $\mathbb{U}_{n}$ is defined by induction as follows:%
\begin{eqnarray*}
\mathbb{U}_{0} &=&\mathbb{R}; \\
\mathbb{U}_{n+1} &=&\mathbb{U}_{n}\cup \mathcal{P}\left( \mathbb{U}%
_{n}\right)
\end{eqnarray*}
\end{definition}

Here $\mathcal{P}\left( E\right) $ denotes the power set of $E.$ If we
identify the couples with the Kuratowski pairs and the functions and the
relations with their graphs, $\mathbb{U}$ formalizes the intuitive idea of
mathematical universe.\newline
We denote by $\mathcal{F}$ the set of finite subsets of $\mathbb{U}$.
Ordered with the relation of inclusion, $\mathcal{F}$ becomes a direct set;
following the usual nomenclature, we call \textit{net} (with values in $E$)
any function $\varphi :\mathcal{F} \rightarrow E$.\newline
Following \cite{ultra}, we introduce the $\Lambda$-limit axiomatically:

\begin{itemize}
\item \textsf{(}$\Lambda $-\textsf{1)}\ \textbf{Existence Axiom.}\ \textit{%
There is a superreal field} $\mathbb{K}\supset \mathbb{R}$ \textit{such that
for every net }$\varphi :\mathcal{F} \rightarrow \mathbb{R}$\textit{\ there
exists a unique element }$L\in \mathbb{K\ }$\textit{called the}
\textquotedblleft $\Lambda $-limit" \textit{of}\emph{\ }$\varphi .$ \textit{%
The} $\Lambda $-\textit{limit will be denoted by} 
\begin{equation*}
L=\lim_{\lambda \uparrow \mathbb{U}}\varphi (\lambda ).
\end{equation*}%
\textit{Moreover we assume that every}\emph{\ }$\xi \in \mathbb{K}$\textit{\
is the }$\Lambda $-limit\textit{\ of some net}\emph{\ }$\varphi :\mathcal{F}
\rightarrow \mathbb{R}$\emph{. }

\item ($\Lambda $-2)\ \textbf{Real numbers axiom}. \textit{If }$\varphi
(\lambda )$\textit{\ is} \textit{eventually} \textit{constant}, \textit{%
namely} $\exists \lambda _{0}\in \mathcal{F} :\ \forall \lambda \supset
\lambda _{0},\ \varphi (\lambda )=r,$ \textit{then}%
\begin{equation*}
\lim_{\lambda \uparrow \mathbb{U}}\varphi (\lambda )=r.
\end{equation*}

\item ($\Lambda $-3)\ \textbf{Sum and product Axiom}.\ \textit{For all }$%
\varphi ,\psi :\mathcal{F} \rightarrow \mathbb{R}$\emph{: }%
\begin{eqnarray*}
\lim_{\lambda \uparrow \mathbb{U}}\varphi (\lambda )+\lim_{\lambda \uparrow 
\mathbb{U}}\psi (\lambda ) &=&\lim_{\lambda \uparrow \mathbb{U}}\left(
\varphi (\lambda )+\psi (\lambda )\right); \\
\lim_{\lambda \uparrow \mathbb{U}}\varphi (\lambda )\cdot \lim_{\lambda
\uparrow \mathbb{U}}\psi (\lambda ) &=&\lim_{\lambda \uparrow \mathbb{U}%
}\left( \varphi (\lambda )\cdot \psi (\lambda )\right).
\end{eqnarray*}
\end{itemize}

\begin{theorem}
The axioms ($\Lambda $-1)\textsf{,}($\Lambda $-2),($\Lambda $-3) are
consistent.
\end{theorem}

\textbf{Proof}. This is the content of Theorem 7 in \cite{ultra}.

$\square $

\bigskip

We say that a net $\varphi :\mathcal{F}\rightarrow \mathbb{U}$ is \textbf{%
bounded} if

\begin{center}
$\exists n\in \mathbb{N}$ such that, $\forall \lambda \in \mathcal{F},$ $%
\varphi (\lambda )\in \mathbb{U}_{n}$.
\end{center}

The notion of $\Lambda $-limit can be extended to bounded nets by induction
on $n$: for $n=0,$ $\lim\limits_{\lambda \uparrow \mathbb{U}}\varphi
(\lambda )$ is defined by the axioms \textsf{(}$\Lambda $-\textsf{1),}($%
\Lambda $-2),($\Lambda $-3); so by induction we may assume that the limit is
defined for $n-1$ and we define it for a net $\varphi :\mathcal{F}%
\rightarrow \mathbb{U}_{n}$ as follows:%
\begin{equation*}
\lim_{\lambda \uparrow \mathbb{U}}\varphi (\lambda )=\left\{ \lim_{\lambda
\uparrow \mathbb{U}}\psi (\lambda )\ |\ \psi :\mathcal{F}\rightarrow \mathbb{%
U}_{n-1}\ \text{and},\ \forall \lambda \in \Lambda ,\ \psi (\lambda )\in
\varphi (\lambda )\right\}
\end{equation*}%
A set that is a $\Lambda $-limit of sets is called \textbf{internal}. The $%
\Lambda $-limit provides a way to extend subset of $\mathbb{R}$ and
functions defined on (subsets of) $\mathbb{R}$ to $\mathbb{K}$:

\begin{definition}
Given a set $E\subset \mathbb{R}$ let $c_{E}:\mathcal{F}\rightarrow \mathbb{U%
}$ be the net such that $\forall \lambda \in \mathcal{F}$ $c_{E}(\lambda )=E 
$. Then%
\begin{equation*}
E^{\ast }:=\lim_{\lambda \uparrow \mathbb{U}}c_{E}(\lambda )=\ \left\{
\lim_{\lambda \uparrow \mathbb{U}}\psi (\lambda )\ |\ \psi (\lambda )\in
E\right\}
\end{equation*}%
is called \textbf{natural extension }of $E$.
\end{definition}

Using the above definition we have that 
\begin{equation*}
\mathbb{K}=\mathbb{R}^{\ast }.
\end{equation*}%
A function $f$ can be extended by identifying $f$ and its graph, and this
extension satisfies the following properties:

\begin{theorem}
For every sets $A,B\in\mathbb{U}$, the \textbf{natural extension} of a
function%
\begin{equation*}
f:A\rightarrow B
\end{equation*}%
is a function 
\begin{equation*}
f^{\ast }:A^{\ast }\rightarrow B^{\ast };
\end{equation*}%
moreover for every $\varphi :\Lambda \cap \mathcal{P}\left( A\right)
\rightarrow A,$ we have that%
\begin{equation*}
\lim_{\lambda \uparrow \mathbb{U}}\ f(\varphi (\lambda ))=f^{\ast }\left(
\lim_{\lambda \uparrow \mathbb{U}}\varphi (\lambda )\right) .
\end{equation*}
\end{theorem}

A property that is natural to ask for the $\Lambda $-limit of a net $\varphi 
$ is that some properties of the limit can be deduced from the properties of 
$\varphi $. This is ensured by the following important theorem:

\begin{theorem}
\label{limit}\textbf{(Leibnitz Principle)} Let $\mathcal{R}$ be a relation
in $\mathbb{U}_{n}$ for some $n\geq 0$ and let $\varphi $,$\psi:\mathcal{\ F}%
\rightarrow\mathbb{U}_{n}$. If 
\begin{equation*}
\forall \lambda \in \mathcal{F},\ \varphi (\lambda )\mathcal{R}\psi (\lambda
)
\end{equation*}%
then%
\begin{equation*}
\left( \underset{\lambda \uparrow \mathbb{U}}{\lim }\varphi (\lambda
)\right) \mathcal{R}^{\ast }\left( \underset{\lambda \uparrow \mathbb{U}}{%
\lim }\psi (\lambda )\right)
\end{equation*}
\end{theorem}

\bigskip

The last key concept that we need is that of hyperfinite set:

\begin{definition}
An internal set is called \textbf{hyperfinite} if it is the $\Lambda $-limit
of finite sets.
\end{definition}

All the internal finite sets are hyperfinite, but there are hyperfinite sets
which are not finite, e.g. the set%
\begin{equation*}
\mathbb{R}^{\circ }:=\ \underset{\lambda \uparrow \mathbb{U}}{\lim }(\mathbb{%
R}\cap \lambda )
\end{equation*}%
is not finite. The hyperfinite sets are very important since, by Leibnitz
Principle, they inherit many properties of finite sets; e.g., $\mathbb{R}%
^{\circ }$ has a maximum and a minimum element, and every internal function
(i.e. a function such that its graph is an internal set) 
\begin{equation*}
f:\mathbb{R}^{\circ }\rightarrow \mathbb{R}^{\ast }
\end{equation*}%
has a maximum and a minimum as well. Intuitively, hyperfinite sets can be
thought as having an hyperfinite number $\beta $ of elements, where $\beta $
is an element of $\mathbb{N}^{\ast }$.\newline
Given a set $A\in \mathbb{U}$ we denote by $A^{\circ }$ its hyperfinite
extension:

\begin{equation*}
A^{\circ}=\lim_{\lambda\uparrow\mathbb{U}} (\lambda\cap A).
\end{equation*}

By this construction, if a hyperfinite set consists of numbers, or vectors,
it is possible to add \textbf{all} its elements. Let%
\begin{equation*}
A:=\ \underset{\lambda \uparrow \mathbb{U}}{\lim }A_{\lambda }
\end{equation*}%
be a hyperfinite set; the hyperfinite sum of the elements of $A$ is defined
as follows: 
\begin{equation*}
\sum_{a\in A}a=\ \underset{\lambda \uparrow \mathbb{U}}{\lim }\sum_{a\in
A_{\lambda }}a
\end{equation*}%
In particular, if $A=\{a_{1},...,a_{\beta }\}$ consists of $\beta $
elements, with $\beta \in \mathbb{N}^{\ast }$, we use the notation 
\begin{equation*}
\sum_{a\in A}a=\sum_{j=1}^{\beta }a_{j}.
\end{equation*}

\subsection{Definition of the ultrafunctions}

Let $\Omega $ be a subset of $\mathbb{R}^{N}$, and let $V_{G}\left( \Omega
\right) \ $ be a vector space such that $\mathcal{D}(\overline{\Omega }%
)\subseteq V_{G}\left( \Omega \right) \subseteq \mathcal{C}(\overline{\Omega 
})\cap L^{2}(\Omega )$.\newline
Let $\varphi _{V_{G}(\Omega )}$ be the net such that, for every $\lambda \in 
\mathcal{F}$, $\varphi _{V_{G}(\Omega )}(\lambda )=V_{\lambda }(\Omega )$,
where

\begin{center}
$V_{\lambda }(\Omega )=Span(V_{G}(\Omega )\cap \lambda ).$
\end{center}

\begin{definition}
The \textbf{set of ultrafunctions generated by $V_{G}(\Omega )$} is%
\begin{equation*}
V(\Omega )=\lim\limits_{\lambda \uparrow \mathbb{U}}V_{\lambda }(\Omega
)=Span(V_{G}(\Omega )^{\circ });
\end{equation*}%
any element $u(x)$ of $V(\Omega )$ is called \textbf{ultrafunction} and $%
V_{G}\left( \Omega \right) $ is called the generating space.
\end{definition}

Observe that, being the $\Lambda$-limit of a net of vectorial spaces of
finite dimensions, $V(\Omega)$ is a vectorial space of hyperfinite
dimension. Its dimension, that we denote by $\beta$, is

\begin{center}
$\beta=\lim\limits_{\lambda\uparrow\mathbb{U}} dim(V_{\lambda}(\Omega))$.
\end{center}

The ultrafunctions are $\Lambda $-limits of continuous functions in $%
V_{\lambda }(\Omega )$, so they are internal functions 
\begin{equation*}
u:\Omega ^{\ast }\rightarrow \mathbb{C}^{\ast }.
\end{equation*}%
(we recall that a function is called "internal" if it is a $\Lambda $-limit
of functions).

Notice that $V(\Omega )$ inherits an Euclidean structure that is the $%
\Lambda $-limit of the Euclidean structure of every space $V_{\lambda
}(\Omega )$ given by the usual $L^{2}\left( \Omega \right) $ scalar product;
also, since $V(\Omega )$ is a subset of $L^{2}(\Omega )^{\ast }$, it can be
equipped with the following scalar product%
\begin{equation*}
\left( u,v\right) =\int_{\Omega }^{\ast }u(x)\overline{v(x)}\ dx.
\end{equation*}%
where $\int_{\Omega }^{\ast }$ is the natural extension of the Lebesgue
integral considered as a functional.\newline
Being a vectorial space of hyperfinite dimension, $V(\Omega )$ admits an
hyperfinite orthonormal basis $\{e_{i}(x)\mid i\leq \beta \}$. Having fixed
a basis, we can make two important constructions in an explicit form. The
first is the extension to $V(\Omega )$ of continuous functions $f(x)$ such
that

\begin{equation*}
\forall v(x)\in V(\Omega )\ \ -\infty <\int_{\Omega }^{\ast }f^{\ast
}(x)v(x)dx<+\infty .
\end{equation*}

Let $f(x)$ be such a function, and let $\Phi$ denote the orthogonal
projection

\begin{center}
$\Phi: \mathcal{C}(\Omega)^{\ast}\rightarrow V(\Omega)$.
\end{center}

We call \textbf{canonical extension} of $f(x)$ the ultrafunction

\begin{center}
$f_{\Phi}(x)=\Phi(f^{\ast}(x))$.
\end{center}

Observe that $f_{\Phi}=f^{\ast}\Leftrightarrow f(x)\in V_{G}(\Omega)$ as
expected, and that for every function $f(x)$ the following holds:

\begin{equation*}
\forall v(x)\in V(\Omega), \int^{\ast}
f^{\ast}(x)v(x)dx=\int^{\ast}f_{\Phi}(x)v(x)dx.
\end{equation*}

In terms of the basis $\{e_{i}(x)\mid i\leq\beta\}$, the operator $\Phi$ has
the following expression:

\begin{equation}
\Phi(f(x))=f_{\Phi}(x)=\sum\limits_{i=1}^{\beta} \left(\int^{\ast}
f^{\ast}(\xi)e_{i}(\xi)d\xi\right) e_{i}(x).
\end{equation}

The second important construction regards the Dirac delta functions:

\begin{theorem}
Given a point $q\in \Omega ,$ there exists a unique function $\delta _{q}(x)$
in $V(\Omega )$ such that%
\begin{equation}
\forall v\in V(\Omega ),\ \int^{\ast }\delta _{q}(x)v(x)\ dx=v(q).
\label{ddel}
\end{equation}
\end{theorem}

\textbf{Proof.} The proof can be found in \cite{ultra}, Theorem 23.\newline
$\square$\newline

$\delta _{q}(x)$ is called the \textbf{Dirac ultrafunction} in $V(\Omega )$
concentrated in $q.$ In terms of the basis $\{e_{i}(x)\mid i\leq\beta\}$,
the $\delta_{q}$ has the following expression:

\begin{equation}
\delta_{q}(x)=\sum_{i=1}^{\beta} e_{i}(q) e_{i}(x),
\end{equation}

which validity can be checked with a direct calculation.

\begin{remark}
We observe that, in the context of ultrafunctions, the Dirac ultrafunctions
are actual functions, while in the classical theory of functions they are
distributions. For example, in the ultrafunction context it makes perfect
sense to consider objects like $\delta_{q}(x)^{2}$, $\delta_{q}(x)-1$, $%
\delta_{q}(x)\cdot\delta_{q^{\prime}}(x)$ and so on.
\end{remark}

\section{The model problem}

In this section we want to solve the problem described in the introduction
via a "natural" approach that can not be applied in the classical framework,
while it can be applied in the ultrafunction setting. We begin by describing
the Dirichlet problems in the framework of ultrafunctions.

\subsection{The Dirichlet problem}

Let $\Omega$ be a bounded open set in $\mathbb{R}^{N}$, and consider the
Dirichlet problem:%
\begin{equation}  \label{1}
\left\{ 
\begin{array}{cc}
u\in \mathcal{C}^{2}_{0}(\Omega) &  \\ 
-\Delta u=f(x) & \text{for}\ \ x\in \Omega%
\end{array}%
\right.
\end{equation}

When $\partial \Omega $ and $f(x)$ are smooth problem (\ref{1}) has a unique
solution. Otherwise, in the classical Sobolev approach, problem (\ref{1}) is
trasformed in the following:

\begin{equation}  \label{2}
\left\{ 
\begin{array}{c}
u\in H_{0}^{1}(\Omega ) \\ 
-\Delta u=f(x)%
\end{array}%
\right.
\end{equation}

Problem (\ref{2}) has a unique solution whenever $\Omega $ is a bounded open
set and $f(x)$ is in $H^{-1}(\Omega )$; in this case the equation $-\Delta
u=f$ is required to be satisfied in a weak sense:%
\begin{equation*}
-\int_{\Omega }u\Delta \varphi \ dx=\int_{\Omega }f\varphi \ dx\ \ \forall
\varphi \in \mathcal{D}(\Omega )
\end{equation*}

Also, the solution $u(x)$ given by this procedure is not a function but an
equivalence class of functions defined a.e. in $\Omega .$\newline
In the approach with ultrafunctions let $V_{0}^{2}(\Omega )$ be the space of
ultrafunctions generated by $C_{0}^{2}(\Omega )$. Problem (\ref{1}) can be
rewritten as follows: 
\begin{equation}
\left\{ 
\begin{array}{cc}
u\in V_{0}^{2}(\Omega ) &  \\ 
-\Delta _{\Phi }u=f(x) & \text{for}\ \ x\in \text{\ }\Omega ^{\ast }%
\end{array}%
\right.  \label{3}
\end{equation}%
where $\Delta _{\Phi }=\Phi \circ \Delta ^{\ast }:V_{0}^{2}(\Omega
)\rightarrow V_{0}^{2}(\Omega )$.\newline
Observe that now we are solving the problem in an hyperfinite space, and by
Leibnitz Principle it follows that there is an unique solution for every $%
f(x)\in V_{0}^{2}(\Omega )$ (for the details see \cite{ultra}, Theorem 27).
The idea of the proof is the following. The solution can be constructed by
first finding a solution $u_{\lambda }(x)$ in each finite dimensional space $%
(V_{0}^{2}(\Omega ))_{\lambda }=Span(C_{0}^{2}(\Omega )\cap \lambda )$, and
then taking the $\Lambda $-limit

\begin{equation*}
\overline{u}(x)=\lim\limits_{\lambda \uparrow \mathbb{U}}u_{\lambda }(x).
\end{equation*}

The solution $\overline{u}(x)$ is an ultrafunction defined for every $x\in
\Omega ^{\ast }$ and, since%
\begin{equation*}
\forall x\in \partial \Omega ,\ \forall \lambda \in \mathcal{F}\cap
C_{0}^{2}(\Omega ),\ u_{\lambda }(x)=0,
\end{equation*}%
it follows by Leibniz principle that $\forall x\in \partial (\Omega ^{\ast
}),\ \overline{u}(x)=0.$

So $\overline{u}(x)$ satisfies the pointwise boundary condition, a result
that is not true in the Sobolev approach. Finally, when problem (\ref{1})
has a solution $s(x)\in \mathcal{C}^{2}(\overline{\Omega }),$ then 
\begin{equation*}
\overline{u}(x)=s^{\ast }(x)
\end{equation*}%
and, when problem (\ref{2}) has a solution $g(x)\in H_{0}^{1}(\Omega ),$
then we have that 
\begin{equation*}
\int_{\Omega }g(x)v(x)\ dx\sim \int_{\Omega }^{\ast }\bar{u}(x)v(x)\ dx\ \
\forall v(x)\in \mathcal{C}_{0}^{2}(\overline{\Omega })
\end{equation*}

\subsection{A solution by mean of ultrafunction}

Now let us consider a minimization problem inspired by the one which we have
discussed in the introduction. Let $\Omega \subseteq \mathbb{R}^{N}$ be an
open bounded set; we want to find a function $u$ defined in $\Omega $ (with $%
u=0\ $on $\partial \Omega $) and a point $q\in \overline{\Omega }$ which
minimize the functional%
\begin{equation*}
E(u,q)=\frac{1}{2}\int_{\Omega }|\nabla u(x)|^{2}dx+u(q)
\end{equation*}%
It is well known that this problem has no solution in $\mathcal{C}%
_{0}^{2}\left( \overline{\Omega }\right) $ and it makes no sense in the
space of distributions. On the contrary it is well defined and it has a
solution in $V_{0}^{2}(\Omega ).$ More exactly, we have the following result:

\begin{theorem}
For every point $q\in \overline{\Omega }^{\ast },$ the Dirichlet problem%
\begin{equation*}
\left\{ 
\begin{array}{cc}
\Delta _{\Phi }u=\delta _{q} & \text{for}\ \ x\in \Omega ^{\ast } \\ 
u(x)=0 & \text{for\ }x\in \partial \Omega ^{\ast }%
\end{array}%
\right.
\end{equation*}%
has a unique solution $u_{q}\in V_{0}^{2}(\Omega )$ whose energy $%
E(u_{q},q)\in \mathbb{R}^{\ast }$ is an infinite number; moreover there
exists $q_{0}\in \Omega ^{\ast }$ such that 
\begin{equation*}
E(u_{q_{0}},q_{0})=\ \underset{q\in \overline{\Omega }^{\ast }}{\min }%
E(u_{q},q)=\min\limits_{q\in \overline{\Omega }^{\ast },\ u\in
V_{0}^{2}(\Omega )}E(u,q).
\end{equation*}
\end{theorem}

\textbf{Proof:} First of all we observe that%
\begin{equation*}
\min\limits_{q\in \overline{\Omega }^{\ast }}E(u_{q},q)=\min\limits_{q\in 
\overline{\Omega }^{\ast },\ u\in V_{0}^{2}(\Omega )}E(u,q)
\end{equation*}%
since every stationary point $(\overline{u},\overline{q})$ of $E(u,q)$
satisfies $-\Delta _{\Phi }(\overline{u})=\delta _{\overline{q}}$.

To minimize $E(u,q)$ we use the Feado-Galerkin method, namely the finite
dimensional reduction. First of all, for every $\lambda \in \mathcal{F}$, we
solve the following problem in $(V_{0}^{2}(\Omega ))_{\lambda
}=Span(C_{0}^{2}(\Omega )\cap \lambda )$:

\begin{equation}
\left\{ 
\begin{array}{cc}
u\in (V_{0}^{2}(\Omega ))_{\lambda } &  \\ 
\int \Delta u\ v\ dx=\int \delta _{q}v\ dx & \text{for every}\ \ v\in
(V_{0}^{2}(\Omega ))_{\lambda }%
\end{array}%
\right.  \label{sol}
\end{equation}%
This problem has a unique solution $u_{q,\lambda }(x)$ for every $\lambda
\in \mathcal{F}\cap \mathcal{C}_{0}^{2}(\Omega )$, since $(V_{0}^{2}(\Omega
))_{\lambda }$ is a nonempty finite-dimensional vectorial space. We show
that this solution depends continuosly on $q$. Consider the linear operator 
\begin{equation*}
-\Delta _{\lambda }:(V_{0}^{2}(\Omega ))_{\lambda }\rightarrow
(V_{0}^{2}(\Omega ))_{\lambda }
\end{equation*}%
that associate to every $u$ of $(V_{0}^{2}(\Omega ))_{\lambda }$ the unique
element $-\Delta _{\lambda }(u)$ such that:

\begin{equation}
\forall v\in (V_{0}^{2}(\Omega ))_{\lambda },\ \int_{\Omega }-\Delta
_{\lambda }u\ v\ dx=\int_{\Omega }-\Delta u\ v\ dx.  \label{proiez}
\end{equation}

So, $-\Delta _{\lambda }u$ is the orthogonal projection of $-\Delta u$ on $%
(V_{0}^{2}(\Omega ))_{\lambda }$. Observe that $-\Delta _{\lambda }$ is a
linear operator that acts on a finite dimensional vector space with $%
Ker(-\Delta _{\lambda })=\{0\}$, so it is invertible.\newline
Now, let $e_{1}(x),...,e_{n}(x)$ be an orthogonal base of $(V_{0}^{2}(\Omega
))_{\lambda }$, and consider the function $k:\overline{\Omega }\rightarrow
(V_{0}^{2}(\Omega ))_{\lambda }$ that associates to every point $q\in 
\overline{\Omega }$ the unique function $\delta _{q,\lambda }\in
(V_{0}^{2}(\Omega ))_{\lambda }$ defined as follows:

\begin{equation*}
\delta _{q,\lambda }(x)=\sum_{i=1}^{n}e_{i}(q)e_{i}(x).
\end{equation*}%
Observe that, by definition, $\forall v \in (V_{0}^{2}(\Omega ))_{\lambda }$
we have 
\begin{equation*}
\int_{\Omega }\delta _{q,\lambda }v\ dx=\int_{\Omega
}\sum_{i=1}^{n}e_{i}(q)e_{i}(x)v(x)=\sum_{i=1}^{n}e_{i}(q)\int_{\Omega
}e_{i}(x)v(x)=v(q),
\end{equation*}%
and, since $v(q)=\int_{\Omega }\delta _{q}(x)v(x)\ dx$, we have

\begin{equation}
\forall v\in (V_{0}^{2}(\Omega ))_{\lambda },\ \int_{\Omega }\delta
_{q,\lambda }(x)v(x)\ dx=\int_{\Omega }^{\ast }\delta _{q}(x)v^{\ast }(x).
\label{deltaproiez}
\end{equation}%
Let $u_{q,\lambda }(x)$ be a solution to \ref{sol}. Then, since $-\Delta
_{\lambda }$ is invertible and (\ref{proiez}) and (\ref{deltaproiez}) hold,
we have

\begin{equation*}
u_{q,\lambda }(x)=\Delta _{\lambda }^{-1}\circ k(q).
\end{equation*}%
Since, as observed, $k$ and $(-\Delta _{\lambda })^{-1}$ are continuous
functions, it follows that $u_{q,\lambda }$ depends continuosly on $q$. Thus
also 
\begin{equation*}
F_{\lambda }(q)=E_{\lambda }(u_{q},q)=\frac{1}{2}\int_{\Omega }|\nabla
u_{q,\lambda }(x)|^{2}dx+u(q)
\end{equation*}%
is continuous.\newline
Since $\overline{\Omega }$ is compact, $F_{\lambda }(q)$ has a minimizer
which we denote by $q_{\lambda }.$\newline
Now let

\begin{equation*}
\overline{q}=\lim_{\lambda \uparrow \mathbb{U}}q_{\lambda }
\end{equation*}%
and%
\begin{equation*}
u_{\overline{q}}=\lim_{\lambda \uparrow \mathbb{U}}u_{q_{\lambda },\lambda }.
\end{equation*}%
By Leibniz Principle, $(u_{\overline{q}},\overline{q})$ is the minimizer of $%
E(u,q)$ in $\overline{\Omega }^{\ast }$. Let us see that $\overline{q}\in
\Omega ^{\ast }.$ By definition of Dirac ultrafunction we have that, for all 
$q\in \partial \Omega ,$ $\delta _{q}=0$, so $u_{q}(x)=0$ and $E(u_{q},q)=0$%
, while $E(u_{q},q)<0$ for every $q\in \Omega ^{\ast }$. So $\overline{q}\in
\Omega ^{\ast }$.\newline
$\square $\bigskip

\begin{remark}
A similar problem that can be studied with the same technique is the problem
of a electrically charged pointwise free particle in a box. Representing the
box with an open bounded set $\Omega \subseteq \mathbb{R}^{3}$, denoting by $%
u_{q}$ the electrical potential generated by the particle placed in $q\in
\Omega $, then $u_{q}$ satisfies the Dirichlet problem

\begin{equation*}
\left\{ 
\begin{array}{cc}
u\in \mathcal{C}_{0}^{2}(\Omega ) &  \\ 
\Delta _{\Phi }u=\delta _{q} & \text{for}\ \ x\in \Omega .%
\end{array}%
\right.
\end{equation*}%
The equilibrium point would be the point $q_{0}\in \overline{\Omega }^{\ast
} $ that minimizes the electrostatic energy which is given by%
\begin{equation*}
E_{el}(q)=\frac{1}{2}\int_{\Omega }|\nabla u_{q}(x)|^{2}dx.
\end{equation*}%
Notice that, 
\begin{equation*}
E_{el}(q)=\int_{\Omega }\delta _{q}(x)u_{q}(x)dx-\frac{1}{2}\int_{\Omega
}|\nabla u_{q}(x)|^{2}dx,
\end{equation*}%
namely, on the solution, the electrostatic energy is the opposite than the
energy of a membrane-like problem in $\mathbb{R}^{3}$. In order to solve
this problem we notice that, by definition of Dirac ultrafunction (\ref{ddel}%
), we have that, for all $q\in \partial \Omega ,$ $\delta _{q}=0.$ So $%
E_{el}(q)\geq 0$ and $E_{el}(q)=0$ if and only if $q\in \partial \Omega .$
More precisely we have that

\begin{itemize}
\item $E_{el}(q)\ $is infinite if the distance between $q\ $and $\partial
\Omega $ is larger than some positive real number;

\item $E_{el}(q)$ is positive but not infinite for some $q$ infinitely close
to $\partial \Omega ;$

\item $E_{el}(q)=0$ if and only if $q\in \partial \Omega .$
\end{itemize}
\end{remark}

\bigskip

\end{document}